# IDENTIFICATION OF TARGET SYSTEM OPERATIONS. DETERMINATION OF THE VALUE OF THE COMPLEX COSTS OF THE TARGET OPERATION

*Запропоновано новий підхід до визначення показника «комплексні витрати (ресурсомісткість)». Показано, що даний показник має ряд переваг по відношенню до витрат операції. Встановлено, що тенденція зміни величини ресурсоємності операції при зміні витрат енергетичного продукту, тривалості зв'язаного стану інших технологічних продуктів і зміні величини цільового продукту операції не суперечить загальній концепції втрат управління*

*Ключові слова: комплексні витрати операції, ресурсомісткість операції, витрати енергетичного продукту, витрати операції*

*Предложен новый подход к определению показателя «комплексные затраты (ресурсоёмкость)». Показано, что данный показатель имеет ряд преимуществ по отношению к затратам операции. Установлено, что тенденция изменения величины ресурсоёмкости операции при изменении расхода энергетического продукта, продолжительности связанного состояния других технологических продуктов и изменении величины целевого продукта операции не противоречит общей концепции потерь управления*

*Ключевые слова: комплексные затраты операции, ресурсоёмкость операции, расход энергетического продукта, потери операции*

**I. Lutsenko**
PhD, Professor
Department of Electronic Devices
Kremenchuk Mykhailo Ostrohradshyi
National University
Pervomaiskaya str., 20,
Kremenchuk, Ukraine, 39600
E-mail: delo-do@i.ua

## 1. Introduction

Determination of the time of the actual completion of the target operation [1] outlines the time range of the study, in which data processing of the operation is carried out. This processing is performed in order to obtain parameters that allow to identify the investigated operation.

One of such basic traditional indicators are costs of operation. The concept of "costs" is rather a tribute to the influence of the accounting terms. This concept, for lack of a suitable alternative, is widely used not only by economists, but also specialists in technical systems, as well as cyberneticists. For example, consumption, the need for raw materials and costs of operation are generally considered synonymous.

The difference between these concepts can be evaluated if we turn to the field of management. For example, metal scrap recycling operation is accompanied by energy product consumption and requires processing a well-defined volume of metal scrap. Energy product consumption depends on the energy product feed intensity and is a management function.

However, when it comes to management and optimization, a popular approach is the minimization of energy costs. The question that the minimization of the energy product consumption leads to a change in operation time and, as a rule, to its significant delaying remains «behind the scenes». At the same time, delaying of the operation is associated with the unproductive use of the raw product, and the struggle for minimum energy costs does not take into account this fact.

In this sense, the task of developing a fundamentally different, cybernetic indicator, which allows to take into account the influence of the change in the energy product consumption, duration of binding raw products of the operation on the operation time and effectiveness is interesting. Since the operation is carried out in order to obtain the result, but not with the purpose of saving energy products.

## 2. Analysis of the literature data and problem statement

Cost reduction was and is today a very topical issue [2]. Costs are controlled [3], optimized, reduced [4] and saved [5].

At the same time, it is clear that costs or expenses is not an independent indicator of the target operation, since two operations, that have the same costs and different duration, have different efficiency, and costs do not indicate this fact without the use of additional indicators.

A similar problem occurs if two operations with similar costs have different values added (profit). In this case, profitability indicator, based on the costs is introduced, but the time factor does not allow to consider this indicator as an independent pointer to a more efficient operation. More profitable, but durable operation may be less effective in relation to the less profitable and less durable operation.

The above is true for problems with the Pareto optimum [6].

In such circumstances, obtaining practically useful result is usually achieved using mathematical modeling methods in searching for the best solution [7].

### 3. Goal and objectives

The goal is to develop a fundamentally new indicator of the target operation "complex costs (resource intensity)" and obtain an expression for numerical and analytical determination of this new cybernetic category.

For this purpose, the following tasks were solved:
1. Determination of the value of resource consumption at the time of the actual completion of the target operation;
2. Determination of the value of resource return at the time of the actual completion of the target operation;
3. Determination of losses of management in the form of a closed thread of mismatch of resource consumption and resource return thread on the interval from the start of the operation until the actual completion;
4. Determination of the complex losses of management (resource intensity of the operation) in the form of the integral value of losses of management at the time of the actual completion of the target operation.

### 4. Complex costs of the target operation

As was shown [1], any effective target operation begins with the resource consumption process, passing in the resource return process. The resource consumption process can be displayed as a closed thread of tight resources, and resource return process - in the form of the target thread (Fig. 1).

The time of the actual completion of the operation (MFZO) determines the time of compensation of the resource consumption thread by the resource return thread [8]. However, this compensation is carried out only in magnitude. In time, (it can be clearly seen when considering the thread of tight resources $ibe(t)$ and the resulting thread $ide(t)$), these threads are spread over and, consequently, at the MFZO there are irreplaceable system losses of management.

These losses can be defined as the value of the integral function of mismatch at the MFZO. That is, system losses of management as a closed thread of the function of mismatch $dif(t)$ are determined by the difference of the integral function of the thread of tight resources $vbe(t)$ and integral function of the resulting thread $vde(t)$ on the interval from the time of the start of the operation until the actual completion (Fig. 1).

Equivalent approach is to determine the losses of management by mismatching the integral function of the resource consumption thread $vre(t)$ and integral function of the resource return thread $vpe(t)$ on the interval from the start of the operation until the actual completion.

Let us define the function of system losses of management on an example of the study of a simple operation, set in the form of a tuple $\left(RE = -2,\ t_r = 2;\ PE = 4,\ t_p = 8\right)$.

1. We determine the MFZO from the expression

$$t_a = \frac{PE \cdot t_p - |RE| \cdot t_r}{PE - |RE|} = 20 \text{ s}.$$

2. We introduce an auxiliary variable $v$ and denote it as $v \in [t_0; t_a]$.

3. We build a model of the operation on the interval $\left[v_0 = t_0; v_a = t_a\right]$ (Fig. 2).

$$ire(v) = \int_{v_0}^{v} re(v)dv,\quad ipe(v) = \int_{v_0}^{v} pe(v)dv.$$

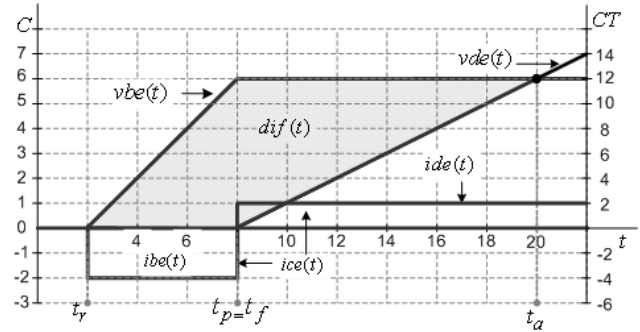

Fig. 1. Determination of system losses of management using the deployed model of the operation ice (t): $t_r$ — the time of the feed of input products of the operation; $t_p$ — the time of the release of output products of the operation; $t_f$ — the time of the physical completion of the target operation; $t_a$ — the time of the actual completion of the target operation

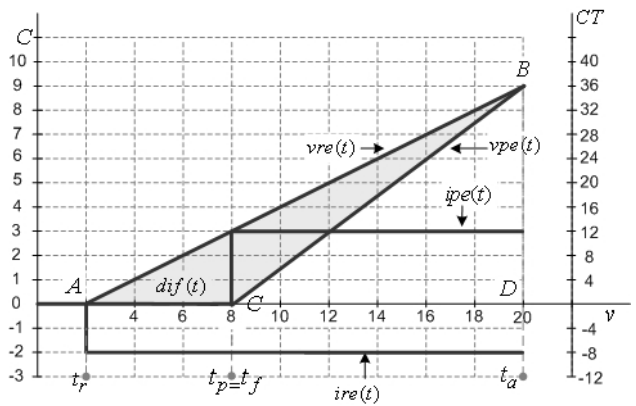

Fig. 2 Determination of system losses of management using the deployed model of the operation in the form of threads ire (t) and ipe(p)

4. We form the integral functions of resource consumption and resource return threads on the interval $\left[v_0 = t_0; v_a = t_a\right]$ (Fig. 2).

$$vre(v) = \int_{v_0}^{v}\left|\int_{v_0}^{v} re(v)dv\right|dv,\quad vpe(v) = \int_{v_0}^{v}\left(\int_{v_0}^{v} pe(v)dv\right)dv.$$

5. We define losses of management as a closed thread of mismatch on the interval from the start of the operation until the actual completion (Fig. 3).

$$dif(v) = \int_{v_0}^{v}\left|\int_{v_0}^{v} re(v)dv\right|dv - \int_{v_0}^{v}\left(\int_{v_0}^{v} pe(v)dv\right)dv.$$

6. We determine the integral function of losses of management (Fig. 4)

$$r(v) = \int_{v_0}^{v} \left[ \int_{v_0}^{v} \left| \int_{v_0}^{v} re(v)dv \right| dv - \int_{v_0}^{v} \left( \int_{v_0}^{v} pe(v)dv \right) dv \right] dv.$$

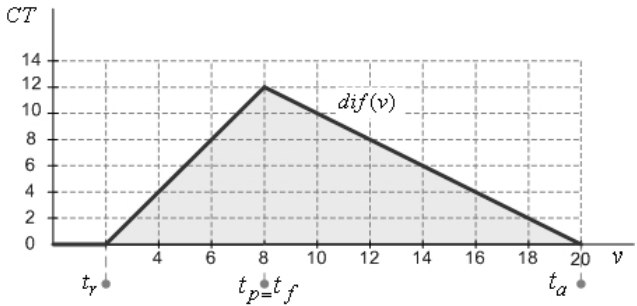

Fig. 3 Uncompensated losses of management in the form of a closed thread dif(v)

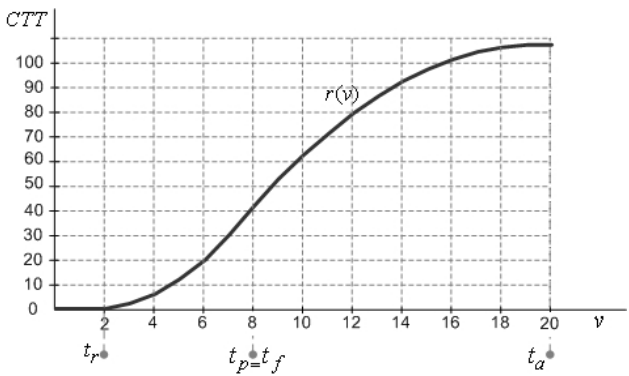

Fig. 4 The function of change of complex losses of management (resource intensity) of the target operation

The value of the resource intensity of the target operation is numerically equal to the value of the function $r(v)$ at the point of MFZO. Therefore, the expression for determining resource intensity will be of the form

$$R = \int_{t_0}^{t_a} \left[ \int_{v_0}^{v} \left| \int_{v_0}^{v} re(v)dv \right| dv - \int_{v_0}^{v} \left( \int_{v_0}^{v} pe(v)dv \right) dv \right] dv, \; v \in [0, t_a]$$

where R is the value of the resource intensity of the operation.

For resource intensity, let us define a unit of measure CTT. Examples of the practical use of formulas for numerical and analytical determination of resource intensity are available on the resource [9].

### 5. Analytical determination of complex costs (resource intensity) of a simple target operation

Since the resource intensity of the simple target operation, based on the geometric interpretation of losses of management, is the area of triangle ABC (Fig. 2), it can be defined as the difference of triangles ABD and CBD.

By substituting scalar values $\text{ire}^*[t_r] \cdot t_r$ and $\text{ipe}^*[t_p] \cdot t_p$ for the functions $\text{ire}^*(t) \cdot t$ and $\text{ipe}^*(t) \cdot t$ we obtain the expressions

$$0 = \text{ire}^*[t_r] \cdot t_r - \text{ire}^*[t_a] \cdot t_a + C,$$
$$0 = \text{ipe}^*[t_p] \cdot t_p - \text{ipe}^*[t_a] \cdot t_a + C.$$

By presenting them in the form of a system of equations and solving with respect to C, we obtain an expression for determining the height BD

$$BD = \frac{\text{ire}^*[t_r] \cdot \text{ipe}^*[t_p] \cdot t_r - \text{ire}^*[t_r] \cdot \text{ipe}^*[t_p] \cdot t_p}{\text{ire}^*[t_r] - \text{ipe}^*[t_p]}.$$

Since we have defined the resource intensity of the simple operation as the difference of the right triangles, the expression for determining it will have the form of

$$R = \frac{1}{2}(t_a - t_r) \cdot BD - \frac{1}{2}(t_a - t_p) \cdot BD.$$

By substituting its value for BD from the expression after the corresponding transformations, we obtain

$$R = \frac{\text{ipe}^*[t_p] \cdot \text{ire}^*[t_r] \cdot (t_r - t_p)^2}{2 \cdot (\text{ipe}^*[t_p] - \text{ire}^*[t_r])}.$$

Given that for simple operations $\text{ipe}^*[t_p]$ is numerically equal to the value PE, and $\text{ire}^*[t_r]$ is numerically equal to the value RE, to determine the numerical value of the resource intensity we can use the expression, which uses the values of the registration signals and moments of their formation

$$R = \frac{PE \cdot |RE| \cdot (t_r - t_p)^2}{2 \cdot (PE - RE)}.$$

For example, for the considered operation (Fig. 1), we obtain

$$R = \frac{PE \cdot |RE| \cdot (t_r - t_p)^2}{2 \cdot (PE - |RE|)} = \frac{3 \cdot 2 \cdot 36}{2} = 108 \; CTT.$$

The value of the resource intensity of the focused operation is numerically equal to the ratio of the product of cost estimates of the module of the value of the input product, the value of the output product and the square of the difference between the time of their registration to the doubled value of the difference between the cost estimates of the module of the value of the input and output products.

As can be seen, the results obtained using a numerical method correspond to the result obtained using the analytical expression for determining resource intensity.

Let us consider the way the resource intensity responds to changes in the parameters of target operations in those cases where the result is guaranteed to be predictable.

Each operation of a set has the cost estimate of input products of the operation (RE – costs), cost estimate of output products of the operation (PE) and the operation time ($T_{op}$). Also, resource intensity (R), which is seen in the diagrams together with costs is calculated for each operation.

In the first set of operations, values RE and PE do not change, and the operation time $T_{op}$ changes. Obviously, the longer the operation time (at fixed RE and PE ), the longer the input products of the operation RE are bound by technological processes and the higher the resource intensity.

Calculation of the resource intensity for the first set of operations (Table 1) confirms this assumption (Fig. 5).

Table 1

Sets of operations to study the relationship of resource intensity with the main traditional indicators of target operations

| Set of operations 1 | | | | Set of operations 2 | | | | Set of operations 3 | | | | |
|---|---|---|---|---|---|---|---|---|---|---|---|---|
| N | RE | PE | T | R | N | RE | PE | T | R | N | RE | PE | T | R |
| 1 | 2 | 3 | 1 | 3 | 1 | 2 | 3 | 3 | 27 | 1 | 2 | 2.5 | 3 | 45 |
| 2 | 2 | 3 | 2 | 12 | 2 | 2.1 | 3 | 3 | 31.5 | 2 | 2 | 2.6 | 3 | 39 |
| 3 | 2 | 3 | 3 | 27 | 3 | 2.2 | 3 | 3 | 37.13 | 3 | 2 | 2.7 | 3 | 34.71 |
| 4 | 2 | 3 | 4 | 48 | 4 | 2.3 | 3 | 3 | 44.36 | 4 | 2 | 2.8 | 3 | 31.5 |
| 5 | 2 | 3 | 5 | 75 | 5 | 2.4 | 3 | 3 | 54 | 5 | 2 | 2.9 | 3 | 29 |
| 6 | 2 | 3 | 6 | 108 | 6 | 2.5 | 3 | 3 | 67.5 | 6 | 2 | 3 | 3 | 27 |
| 7 | 2 | 3 | 7 | 147 | 7 | 2.6 | 3 | 3 | 87.75 | 7 | 2 | 3.1 | 3 | 25.36 |

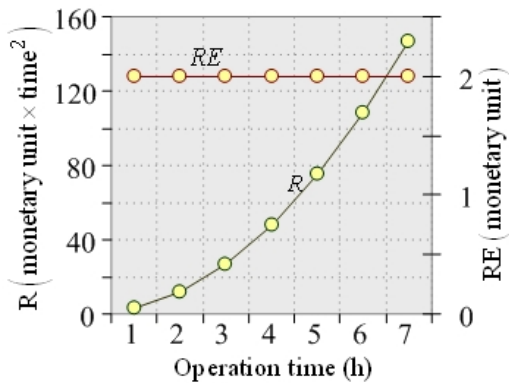

Fig. 5. Change in costs RE and resource intensity R from the target operation time

The second set of operations (Table 1) is characterized by the fact that, from operation to operation, the cost estimate of input products of the operation RE increases and the cost estimate of output products of the operation PE and the operation time $T_{op}$ do not change. In this case, with an increase in the cost estimate of input products of the operation (costs RE ), at constant PE and $T_{op}$ resource intensity of the operation should increase. Calculation of the resource intensity for the second set of operations (see Table) confirms this assumption (Fig. 6).

The third set of operations (Table 1) is characterized by the fact that, from operation to operation, the cost estimate of output products of the operation PE increases and the cost estimate of input products of the operation RE and the operation time $T_{op}$ do not change. Here, the change of the cost estimate of output products of the operation towards an increase speeds up the compensation of the resource consumption of the operation, and resource intensity must, in this case, decrease. Calculation of the resource intensity for the third set of operations (Table 1) confirms this assumption as well (Fig. 7).

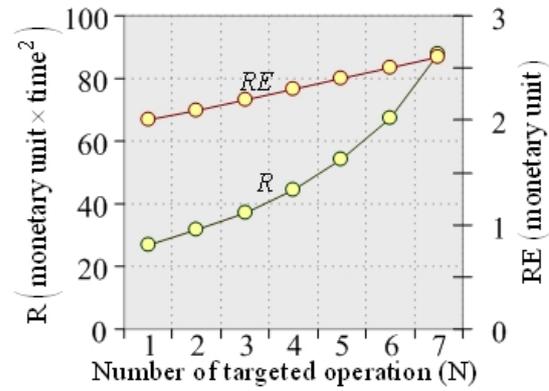

Fig. 6. Change in the resource intensity R from the target operation costs RE

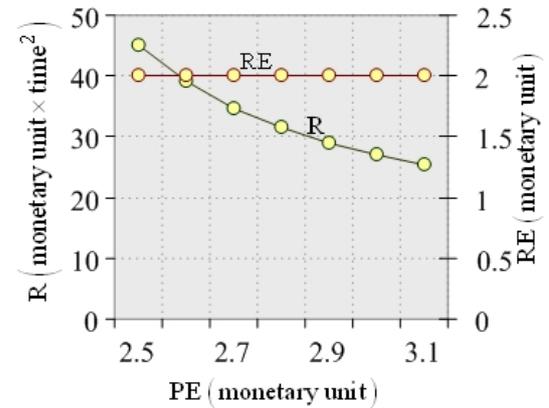

Fig. 7. Change in costs RE and resource intensity R from the expert (cost) estimate of output products PE of the target operation

In conclusion, let us consider the cycle of operations with minimum costs in the process of the proportional change of the operation time (Table 2).

Table 2

Change in the basic parameters of operations in a system having minimum costs

| N | RE | PE | T | R | Prof |
|---|---|---|---|---|---|
| 1 | 2 | 2,5 | 1 | 5 | 575,00 |
| 2 | 1,894 | 2,5 | 1,05 | 4,31 | 663,71 |
| 3 | 1,824 | 2,5 | 1,1 | 4,08 | 706,73 |
| 4 | 1,772 | 2,5 | 1,15 | **4,02** | 728,00 |
| 5 | 1,75 | 2,5 | 1,2 | 4,2 | 718,75 |
| 6 | 1,738 | 2,5 | 1,25 | **4,45** | 701,04 |
| 7 | 1,759 | 2,5 | 1,3 | 5,01 | 655,50 |
| 8 | 1,791 | 2,5 | 1,35 | 5,75 | 603,96 |
| 9 | 1,837 | 2,5 | 1,4 | 6,79 | 544,61 |
| 10 | 1,913 | 2,5 | 1,45 | 8,56 | 465,55 |
| 11 | 2 | 2,5 | 1,5 | 11,25 | 383,33 |

Calculation of resource intensity shows that its minimum does not match the minimum costs (Fig. 8).

Let us analyze the differences among the operations, indicated by the minimum costs (operation N6) and a minimum of resource intensity (operation N4).

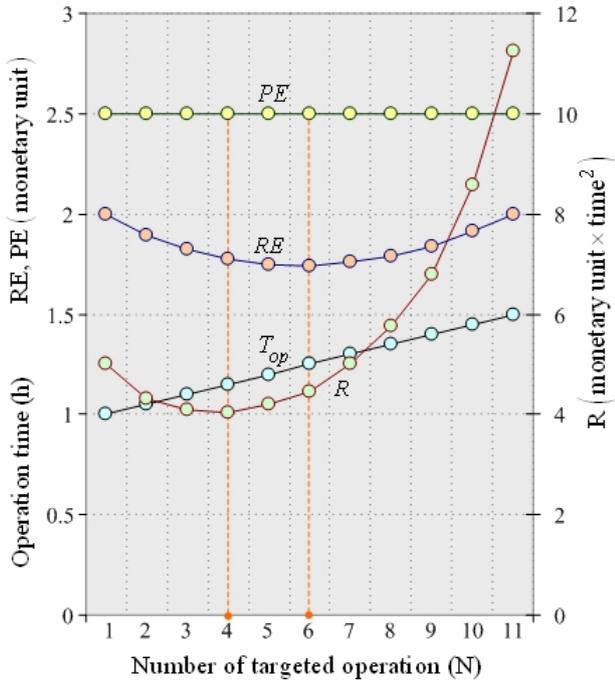

Fig. 8. Resource intensity of the target operation depending on the change in costs and operation time

Thousand N4 type operations, carried out in the cycle lasts for T=Top*1000=1150 hours. During this time, the operation generates the value added (profit) Prof4=1150*(PE-RE)=1150*(2,5-1,772)=728 monetary units.

During the same time, the N6 type operation will be performed I=1150/1.25=920 times. In this case, the N6 type operation generates the value added (profit) Prof6= =920*(PE-RE)=920*(2,5-1,738)=701.04 monetary units.

Thus, the target product Prof4 exceeds the target product Prof6 by 27 monetary units in absolute terms.

Carrying out these calculations for all operations of the set shows that the minimum resource intensity indicates the maximum generated target product (profit) (Fig. 9).

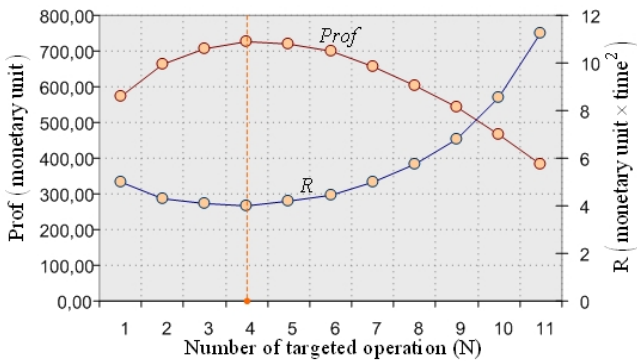

Fig. 9. Change in the accumulated value added (Prof) of operations on the time interval of 1150 hours and resource intensity (R)

As can be seen from the Fig. 9, a change in resource intensity of the target operation is the mirror image of the function of the value added. Consequently, the maximum efficiency of a set of operations corresponds to the operation with minimum resource intensity.

A similar pattern of change in resource intensity of the operation of the management is typical for processes with a batch feed of energy products [10].

## 8. Conclusions

The paper gives a new insight into the resource intensity of the target operation, which is its most important indicator.

It was found that the concept of resource intensity of the target operation is based on such core categories as resource consumption of the operation (consumption of input products of the operation in time in comparable values) and resource return of the operation (generation of output products of the operation in time in comparable values).

Reliance on the functions of resource consumption and resource return, from the start of the target operation until the actual completion has allowed to quantify the value of its resource intensity.

Using the model of simple target operation has allowed to obtain an analytical expression of the resource intensity for the operations, in the study of which distributed parameters of registration functions of input and output products of the operation can be neglected.

Using mathematical modeling methods it was revealed that in the case of a fixed value of expert (cost) estimate of output products of the operation, minimum resource intensity of the target operation indicates the maximum efficiency operation with respect to the target product of the operation.

Using the developed indicator in search optimization systems allows to maximize financial result by increasing the amount of the value added by 5–25 %.

# ДОСЛІДЖЕННЯ БАГАТОПАРАМЕТРИЧНИХ ОБ'ЄКТІВ КОНТРОЛЮ ТА УПРАВЛІННЯ МЕТОДОМ ТРИМІРНОГО ІНТЕГРАЛЬНОГО ФУНКЦІОНАЛУ

*На основі теорії реологічних переходів виконані дослідження технологічних процесів. Показано, що такі процеси можна описати інтегральними імпульсними дельта-функціями Дірака. Приведені аналітичні рівняння, за якими можна розрахувати екстремалі технологічного процесу. Це дозволяє забезпечити максимальну ефективність технологічного процесу при мінімумі енергетичних та матеріальних затрат*

*Ключові слова: технологія, контроль, управління, перенесення, реологія, перехід, дифузія, конвекція, екстремум, оптимізація*

*На основе теории реологических переходов выполнены исследования технологических процессов. Показано, что такие процессы можно описать интегральными импульсными дельта-функциями Дирака. Приведены аналитические уравнения, по которым можно рассчитывать экстремали технологического процесса. Это позволяет обеспечить максимальную эффективность технологического процесса при минимуме энергетических и материальных затрат*

*Ключевые слова: технология, контроль, управление, перенос, реология, переход, диффузия, конвекция, экстремум, оптимизация*

**Й. І. Стенцель**
Доктор технічних наук,
професор, завідувач кафедри*
E-mail: stencel@sti.lg.ua

**О. І. Проказа**
Кандидат технічних наук, доцент*
E-mail: kafKISU.Elena@gmail.com

**К. А. Літвінов**
Аспірант*
E-mail: LitvinovK@yandex.ru

*Кафедра комп'ютерно-інтегрованих систем управління
Східноукраїнський національний університет ім. Володимира Даля
пр. Радянський 59-а,
м. Сєвєродонецьк, Україна, 93400

## 1. Вступ

Технологічні процеси (ТП) хімічних, нафтохімічних, нафтопереробних, теплоенергетичних, фармацевтичних, харчових та інших виробництв відносяться до складних взаємопов'язаних багатопараметричних об'єктів контролю та управління з багатьма вхідними, вихідними та впливовими параметрами. Як правило, такі об'єкти описуються нелінійними диференціальними рівняннями перенесення енергії, маси та кількості руху. У зв'язку з тим, що методів розв'язку нелінійних інтегро-диференціальних рівнянь практично немає, то їх спрощують, приводять до звичайних диференціальних рівнянь, нехтують багатьма фізичними параметрами у деякому діапазоні роботи об'єкта і використовують принцип автономності для основних каналів контролю та управління. При описанні багатопараметричних технологічних об'єктів контролю та управління (ТОКУ), як правило, виходять з рівнянь матеріальних і теплових балансів [1], котрі в подальшому обмежуються тільки їх лінійними складовими. Так як такий метод математичного опису ТОКУ